\documentclass{article}

\title{On polytopes simple in edges}
\author{V. Timorin
\thanks{Partially supported by RFBR 99-01-00245 and CRDF RM1-2086} }
\date{}

\usepackage{amsfonts} \usepackage{amssymb}
\def\R{\mathbb{R}} \def\C{\mathbb{C}} \def\T{\mathbb{T}} \def\Ups{\Upsilon}
\def\d{\partial} \def\Diff{{\rm Diff}} \def\Vol{{\rm Vol}} \def\Oc{{\cal O}}
\def\Mc{{\cal M}} \def\Box{\square} \def\leq{\leqslant} \def\geq{\geqslant}
\def\phi{\varphi} \def\Ker{{\rm Ker}} \def\proof{{\sc Proof. }}
\newtheorem{theorem}{Theorem}[section]
\newtheorem{proposition}[theorem]{Proposition}

\newtheorem{lemma}[theorem]{Lemma}
\newtheorem{conjecture}[theorem]{Conjecture}

\begin{document}
\maketitle

{\small We investigate some combinatorial properties of convex polytopes
simple in edges. For polytopes whose nonsimple vertices are 
located sufficiently far one from another, we prove an analog
of the Hard Lefschetz theorem. It implies Stanley's conjecture for 
such polytopes.}

\section*{Introduction}

In this paper, the word ``polytope'' refers always to a convex polytope.
By {\em $d$-polytope} we mean a polytope in $\R^d$ with nonempty
interior.

Denote by $f_k$ the number of $k$-faces of a $d$-polytope.
It is useful to consider another collection of numbers $\{h_k\}$
that is obtained from $\{f_k\}$ by the linear transformation
$$h_k=\sum_{i\geq k} f_i(-1)^{i-k}{i \choose k},\quad
f_k =\sum_{i\geq k} h_i{i\choose k}.$$
The collection of numbers $f_k$ is called the {\em $f$-vector}, the 
numbers $h_k$ constitute the {\em $h$-vector}.

A vertex of a $d$-polytope is said to be {\em simple} if exactly $d$
facets meet at this vertex. A $d$-polytope is called {\em simple} if all its
vertices are simple. The following relations on the $h$-vector (and hence on
the $f$-vector) hold for a simple polytope
\begin{itemize}
\item{\em Dehn-Sommerville equations} \cite{Gr}:
$h_0=h_d=1$, $h_k=h_{d-k}$,
\item{\em Unimodality condition} \cite{St1,McM}:
$$h_0\leq h_1\leq\cdots\leq h_{[d/2]}\geq h_{[d/2]+1}\geq\cdots
\geq h_{d-1}\geq h_d.$$
\end{itemize}
In this paper, we will study the $h$-vectors of a slightly more general
class of polytopes.

A $d$-polytope is called {\em simple in edges} if each its edge is
incident exactly to $d-1$ facets.
We will prove that for any polytope simple in edges all numbers
$h_{[d/2]}$, $h_{[d/2]+1}$,$\dots$, $h_d$ are nonnegative and
$h_k\leq h_{d-k}$ for $k\leq d/2$.
Polytopes simple in edges appear for
instance as (closures of) fundamental polyhedra of groups generated by
reflections in Lobachevskii spaces. A combinatorial study of
polytopes simple in edges carried out by Khovanskii \cite{Kh1} concluded
the proof of the following important theorem. In a Lobachevskii
space of sufficiently high dimension there are no discrete groups generated by
reflections whose fundamental polyhedron has finite volume.
This statement was inspired by works of Nikulin and Vinberg and reduced to
combinatorics by Prokhorov \cite{Prokh}.
We will give another (more direct) combinatorial proof of Khovanskii's
estimate. 

Suppose that nonsimple vertices of a polytope simple in edges are located
sufficiently far one from another. This means that no facet contains two
nonsimple vertices. Such polytopes will be called {\em polytopes with
infrequent singularities}. Then, as we will show in this paper, the
inequalities
$$h_{[d/2]}\geq h_{[d/2]+1}\geq\cdots\geq h_d$$
still hold. Presumably they are true for any polytope simple in edges.

A polytope is said to be {\em integral} provided all its vertices
belong to the integer lattice. With each integral polytope $\Delta$
one associates the {\em toric variety} $X$ \cite{Kh2, Danilov,
Oda, Fulton}. This is a projective complex algebraic variety, singular
in general. It turns out that the intersection cohomology Betti numbers
of $X$ are combinatorial invariants of
$\Delta$ \cite{F,DL}. Denote them by $Ih_k(\Delta)=\dim IH^k(X,\C)$.
For definition and basic results on intersection cohomology see
\cite{GM1,GM2,BBD}.

For example, for a simple polytope $\Delta$ we have $Ih_k(\Delta)=h_k$
(for simple integral polytopes the associated toric varieties are quasi-smooth
so the intersection cohomology coincides with ordinary cohomology).
From Poincar\'e duality and the Hard Lefschetz theorem in the intersection
cohomology of $X$ it follows that for an integral polytope $\Delta$ 
$$Ih_0(\Delta)=1,\ Ih_k(\Delta)=Ih_{d-k}(\Delta),$$
$$Ih_0(\Delta)\leq Ih_1(\Delta)\leq\cdots\leq Ih_{[d/2]}(\Delta)\geq
Ih_{[d/2]+1}(\Delta)\geq\cdots\geq Ih_d(\Delta).$$
Stanley defined for an arbitrary polytope $\Delta$ a collection of numbers
$Gh_k(\Delta)$ ({\em generalized $h$-vector}) such that
$Gh_k(\Delta)=Ih_k(\Delta)$ for an integral
polytope $\Delta$ \cite{St2}. The definition of $Gh_k$ was motivated
by the calculation of intersection cohomology of toric varieties made
by Bernstein, Khovanskii and MacPherson (independently). 
Stanley proved that $Gh_k=Gh_{d-k}$ and conjectured
that all the inequalities above also hold for arbitrary polytopes
(with $Gh_k$ instead of $Ih_k$). For simple polytopes $Gh_k=h_k$ and
hence Stanley's conjecture is true in this case. We will see that Stanley's
conjecture is
true also for polytopes with infrequent singularities (in this case
$Gh_k=h_k$ for $k\geq d/2$).

In \cite{BL,BBFK} there is a combinatorial description of the intersection
cohomology of toric varieties. This description makes sense for
arbitrary polytopes (not necessarily integral). So for every polytope
$\Delta$ there is the {\em combinatorial intersection cohomology}.
Denote the combinatorial Betti numbers in the same way
$Ih_k(\Delta)$ as in integral case. It is proven in \cite{BL,BBFK} that
$Ih_k=Ih_{d-k}$ and an analog of Poincar\'e pairing is constructed.
Moreover, there is an analog of the Lefschetz operator that coincides
with the ordinary Lefschetz operator in integral case. Presumably a
combinatorial analog of the Hard Lefschetz theorem holds in general. It
would imply that $Ih_k=Gh_k$ so this is a stronger version of Stanley's
conjecture. In this paper, we will prove some variant of the combinatorial
Hard Lefschetz theorem for polytopes with infrequent singularities. 

{\em Acknowledgements.}
I am very grateful to V. A. Lunts for the significant help (he has told me
the statement of theorem \ref{comM}) and to A. G. Khovanskii for useful
discussions.

\section{Cohomology of simple polytopes}

In this section, we recall the geometric definition of the cohomology
of simple polytopes given in \cite{PKh}. 

Let $\Sigma$ be a $d$-polytope. A polytope $\Sigma'$ is said to be
{\em analogous} to $\Sigma$ if there is a one-to-one correspondence between the facets
of $\Sigma$ and $\Sigma'$ such that corresponding facets have the 
the same outer normals and become analogous after being shifted to a common hyperplane. 
By definition any two segments in $\R$ are analogous. Analogous polytopes have the 
same {\em combinatorial type}. This means that faces of analogous polytopes satisfy
the same inclusion-relations. In particular, a polytope analogous to a simple one is
simple. 

Fix any simple $d$-polytope $\Sigma$ in $\R^d$. Consider a polytope
$\Sigma'$ analogous to $\Sigma$. For each facet $\Gamma$ of $\Sigma$ there
is the corresponding (parallel) facet $\Gamma'$ of $\Sigma'$. Let
$\xi_\Gamma$ be a linear functional whose maximal value on $\Sigma$ is
achieved on the facet $\Gamma$. Denote by $H_\Gamma(\Sigma')$ the maximum of
$\xi_\Gamma$ on $\Sigma'$ (of course, this maximum is achieved on
$\Gamma'$). The number $H_\Gamma(\Sigma')$ is called a {\em support number}
of $\Sigma'$.

Move slightly all the facets of $\Sigma$ so that each remains parallel to
itself. Then we get an analogous polytope $\Sigma'$. It follows that we can
vary the support numbers independently (at least while the divergences are
sufficiently small). On the other hand, the support numbers of $\Sigma'$
determine $\Sigma'$. So we can think of the polytope $\Sigma'$ as a function
in the independent parameters $H_\Gamma$.

The volume of $\Sigma'$ turns out to be a  polynomial in $H_\Gamma$. Denote
this polynomial by $\Vol_\Sigma$. Now consider the ring $\Diff$ of all
differential operators with constant coefficients with respect to the support
numbers. Denote by $\d_\Gamma$ the operator of differentiation with respect
to $H_\Gamma$. The ring $\Diff$ is nothing more than the polynomial ring in
the differentiations $\d_\Gamma$. Let $J$ be the ideal in $\Diff$ consisting
of operators $\alpha$ such that $\alpha\Vol_\Sigma=0$. The ideal $J$ is
homogeneous, hence the quotient algebra $A(\Sigma)=\Diff/J$ inherits the
grading. Dimension of the homogeneous component $A^k(\Sigma)$ equals to
$h_k(\Sigma)$ \cite{Tim}.

The ideal $J$ can be described explicitly. It is generated by the
following two groups of differential operators \cite{Tim}:
\begin{itemize}
\item $\d_{\Gamma_1}\cdots\d_{\Gamma_k}$ where $\Gamma_1\cap\cdots\cap
\Gamma_k=\varnothing$,
\item $L_a=\sum H_\Gamma(a)\d_\Gamma$ where $a\in\R^d$ is a point
considered as a limit case of a polytope analogous to $\Sigma$. 
\end{itemize}
The second group is responsible for the translation invariance of volume.

The ring $A(\Sigma)$ models the cohomology ring. If $\Sigma$ is integral,
then $A(\Sigma)$ is indeed isomorphic to the cohomology ring of the
corresponding toric variety. The operator of multiplication by
$L_\Sigma=\sum H_\Gamma(\Sigma)\d_\Gamma$ represents the Lefschetz operator
(that is dual to the hyperplane section opetator in the homology). The
following analog of the Hard Lefschetz theorem holds for simple polytopes.

\begin{theorem}
\label{HLs}
The operator of multiplication by $L_\Sigma^{d-2k}$ establishes an
isomorphism betweeen $A^k(\Sigma)$ and $A^{d-k}(\Sigma)$. 
\end{theorem}

The first proof is due to McMullen \cite{McM}. He uses another description
of cohomology. See also \cite{Tim}. From this theorem it follows that the
$h$-vector of a simple polytope is unimodal, i.e.,
$h_0\leq h_1\leq\cdots\leq h_{[d/2]}$. An element $\alpha\in A^k(\Sigma)$ is
called {\em primitive} if $\alpha L_\Sigma^{d-2k+1}=0$ in $A^k(\Sigma)$ or,
equivalently, the polynomial $\alpha\Vol_\Sigma$ has zero of order $k$ at the
point with coordinates $H_\Gamma(\Sigma)$.
It is easy to see that the space of all order-$k$ primitive elements has
dimension $h_k-h_{k-1}$. The following theorem \cite{McM,Tim} is an analog
of the Hodge--Riemann bilinear relations (and a generalization of
the Brunn--Minkowski inequality):

\begin{theorem}
\label{HR}
For each primitive element $\alpha\in A^k(\Sigma)$
$$(-1)^k\alpha^2L_\Sigma^{d-2k} (\Vol_\Sigma)\geq 0.$$ 
\end{theorem}

Now let us study the relation between $A^k(\Sigma)$ and $A^k(\Gamma)$ where
$\Gamma$ is a facet of $\Sigma$. First note that the polynomial
$\d_\Gamma\Vol_\Sigma$ evaluated at the support numbers of $\Sigma$ gives
the $(d-1)$-volume of $\Gamma$ (it is almost obvious).
On the other hand, the support numbers of $\Gamma$ are certain linear
functions of the support numbers of $\Sigma$ (that can be written down
explicitly, of course). Therefore the polynomial $\Vol_\Gamma$ differs from
$\d_\Gamma\Vol_\Sigma$ by a linear (noninvertible)
change of variables. Given an element $\alpha\in A^k(\Sigma)$ one can find an element 
$\alpha_{(\Gamma)}\in A^k(\Gamma)$ such that $\alpha_{(\Gamma)}\Vol_\Gamma$ differs from 
$\alpha\d_\Gamma\Vol_\Sigma$ by the same change of variables. The map $\alpha\mapsto
\alpha_{(\Gamma)}$ is a surjective homomorphism from $A^k(\Sigma)$ to $A^k(\Gamma)$. 
It is not hard to show that $(L_\Sigma)_{(\Gamma)}=L_\Gamma$.  

\paragraph{An analog of Morse theory.} A {\em general linear function} on
$\Sigma$ is a linear function on $\R^d$ that is nonconstant on any edge of
$\Sigma$. Fix a general linear function $l$. We will view $l$ as a vertical
coordinate and will apply to it the words ``up'' and ``down''. {\em Index of
a vertex} $v$ of $\Sigma$ is the number of edges that go down from $v$.
It is not difficult to prove that the number of vertices of index $k$ 
in $\Sigma$ equals to $h_k$ \cite{Kh1}. In particular, it does not depend on
the choice of $l$.

Let $v$ be a vertex of $\Sigma$. The {\em separatrix} of $v$
is the face of $\Sigma$ spanned by all the edges that go down from $v$.
There is an explicit description of a basis in the cohomology space
$A^k(\Sigma)$ in terms of differential operators \cite{Tim}. Let $F$ be a
face of $\Sigma$. Denote by $\Gamma_1,\dots,\Gamma_k$ all the facets
containing $F$ so $F=\Gamma_1\cap\cdots\cap\Gamma_k$. Let us define the
differential operator
$\d_F=\d_{\Gamma_1}\cdots\d_{\Gamma_k}$ associated to the face $F$.

\begin{theorem}
Fix any general linear function on a simple polytope $\Sigma$.
The operators $\d_F$, where $F$ is a separatrix of $\Sigma$, constitute
a basis in the vector space $A(\Sigma)$. 
\end{theorem}

The decomposition of an element $\alpha\in A(\Sigma)$ with respect to this
basis is called the {\em separatrix decomposition}. Let
$\alpha=\sum a_F\d_F$ be the separatrix decomposition of $\alpha$. A
separatrix $F$ is called the {\em highest separatrix} of $\alpha$ if value
of $\max(l|_F)$ is highest among all the separatrices that enter
the decomposition of $\alpha$ with nonzero coefficients.

\section{Polytopes simple in edges}

Consider a $d$-polytope $\Delta$ simple in edges. Let $v$ be a 
nonsimple vertex of $\Delta$. Let us ``cut off'' the vertex $v$ from
$\Delta$ by a hyperplane $P$ sufficiently close to $v$. More precisely, we
take a hyperplane $P$ such that $v$ lies on one side with respect to $P$ and
all other vertices of $\Delta$ lie on the other side. Denote by $P^+$ the
half-space that does not contain $v$. Then we say that $\Delta\cap P^+$ is
the polytope $\Delta$ with the vertex $v$ {\em cut off}. Now cut off all
the nonsimple vertices (clearly the result does not depend on the order of
our cut-off processes). We get a simple polytope $\Sigma$ which satisfies
the following condition. There is a continuous one-parameter family
$\Sigma_t$ of analogous simple polytopes such that $\Sigma_1=\Sigma$ and
$\Sigma_t\to \Delta$ as $t\to 0$ in the Hausdorff metric. The polytope
$\Sigma$ will be called the {\em standard resolution} of $\Delta$.
Consider a facet of $\Sigma$ that comes from a cutting hyperplane. We call
such a facet an {\em inserted facet}.

Denote by $L_\Delta$ the limit of operators $L_{\Sigma_t}$ in $A(\Sigma)$.
In coordinates, $L_\Delta=\sum H_\Gamma(\Delta)\d_\Gamma$ where
the support numbers $H_\Gamma(\Delta)$ are defined as maximal values
of the functionals $\xi_\Gamma$ restricted to $\Delta$.

Let $\Gamma$ be a facet of $\Sigma$. It is included in a one-parameter
family $\Gamma_t$ of facets of $\Sigma_t$. The limit
$\Gamma_0=\lim_{t\to 0}\Gamma_t$ is a face of $\Delta$. For example,
for an inserted facet we get just a vertex that was cut off.
Note that $(L_\Delta)_{(\Gamma)}$ coincides with $L_{\Gamma_0}$ in
$A(\Gamma)$. In particular, if $\Gamma$ is the inserted facet corresponding
to a vertex $v$ of $\Delta$, then $(L_\Delta)_{(\Gamma)}=L_v=0$ in
$A(\Gamma)$.

\begin{lemma}
\label{hcalc}
For $0<k\leq d$, we have $h_k(\Delta)=h_k(\Sigma)-\sum h_k(\Gamma)$ where
the sum is over all inserted facets $\Gamma$ of $\Sigma$ (we are
assuming that $f_d(\Gamma)=h_d(\Gamma)=0$).
\end{lemma}

\proof This follows from the analogous formula for the $f$-vector:
$f_m(\Delta)=f_m(\Sigma)-\sum f_m(\Gamma)$. $\Box$

\begin{lemma}
\label{indep}
Assume that $\sum\d_\Gamma\alpha_\Gamma=0$ where $\Gamma$ runs over inserted
facets and $\alpha_\Gamma$ are operators of order $k<(d-1)/2$ in $A(\Sigma)$.
Then $\d_\Gamma\alpha_\Gamma=0$ for each $\Gamma$. 
\end{lemma}

\proof Fix a facet $\Gamma'$ of $\Sigma$ and
multiply the equation $\sum\d_\Gamma\alpha_\Gamma=0$ by $\d_{\Gamma'}$.
Observe that the operator $\d_{\Gamma_1}\d_{\Gamma_2}$ is nonzero in
$A(\Sigma)$ if and olny if $\Gamma_1\cap\Gamma_2\ne\varnothing$.
Since the inserted facets are disjoint, we get
$\d_{\Gamma'}^2\alpha_{\Gamma'}=0$.

Denote by $\Delta'$ the polytope $\Delta$ with all nonsimple vertices but
that corresponding to $\Gamma'$ cut off. We are assuming that the
inserted facets of $\Delta'$ are the same as in $\Sigma$. In particular,
$\Sigma$ is a standard resolution of $\Delta'$. Then $L_{\Delta'}-L_\Sigma=
c\d_{\Gamma'}$ where $c$ is a positive number. We know that $L_\Sigma$
descends to a Lefschetz operator on $A^k(\Gamma')$ and $L_{\Delta'}$
descends to zero. Therefore $\d_{\Gamma'}$ represents a negative
Lefschetz operator in $A^k(\Gamma')$.

Thus $L_{\Gamma'}(\alpha_{\Gamma'})_{(\Gamma')}=0$. By the Hard Lefschetz
theorem \ref{HLs} for $\Gamma'$ we conclude that
$(\alpha_{\Gamma'})_{(\Gamma')}=0$, i.e., $\d_{\Gamma'}\alpha_{\Gamma'}=0$.
$\Box$

\begin{lemma}
\label{ineq1}
For $k\leq d/2$, we have $h_k(\Sigma)-\sum h_{k-1}(\Gamma)\geq 0$.  
\end{lemma}

\proof Consider the subspace in $A^k(\Sigma)$ generated by elements
$\d_\Gamma\alpha$ where $\Gamma$ is an inserted facet of $\Sigma$
and $\alpha\in A^{k-1}(\Sigma)$. By lemma \ref{indep} we know that dimension
of this subspace equals to $\sum\dim(\d_\Gamma A^{k-1}(\Sigma))=
\sum\dim(A^{k-1}(\Gamma))=\sum h_{k-1}(\Gamma)$.
Therefore $h_k(\Sigma)-\sum h_{k-1}(\Gamma)$ is nonnegative. $\Box$

\begin{theorem}
\label{ineq2}
The numbers $h_k(\Delta)$ are nonnegative for all $k\geq d/2$.
\end{theorem}

\proof This follows form lemmas \ref{hcalc} and \ref{ineq1} and the
Dehn-Sommerville equations for $\Sigma$ and all $\Gamma$'s. $\Box$

Note that other components of the $h$-vector need not to be nonnegative. For
example, the icosahedron is simple in egdes (as any 3-dimensional polytope)
but $h_1=-7$.

\begin{theorem}
\label{ineq3}
For $k\leq d/2$ we have $h_k(\Delta)\leq h_{d-k}(\Delta)$. 
\end{theorem}

\proof For $k=0$ this is an equality since $h_0=h_d=1$ (it follows from the
Euler theorem). Suppose $k>0$. Then the inequality follows from lemma
\ref{hcalc}, Dehn-Sommerville equations for $\Sigma$ and the unimodality
condition for the inserted facets $\Gamma$. $\Box$

It is known that the Euler theorem $h_0=1$ and the trivial equation $h_d=1$ 
(these equations are true for any convex polytope) are the only linear
relations on the $h$-vector (equivalently, on the $f$-vector) of a polytope
simple in edges \cite{Gr}.
Theorem \ref{ineq2} provides some inequality-type relations.  
Later on we will prove some more subtle inequalities for polytopes
with infrequent singularities.

\paragraph{Some applications.} A. Khovanskii in \cite{Kh1} estimated the 
average number of $k$-dimensional subfaces on a $l$-dimensional face of 
a $d$-polytope simple in edges ($1\leq k<l\leq d/2$). 
Khovanskii's estimate generalized the earlier result
of Nikulin \cite{Nik} (who worked out the case of simple polytopes) and
completed the proof of the following: in Lobachevskii space of dimension
$>995$ there are no discrete groups generated by reflections with
fundamental polyhedron of finite volume.

We will deduce the Khovanskii's estimate from theorems \ref{ineq2} and
\ref{ineq3}.

The following lemma is almost obvious (it can be easily proved by induction):

\begin{lemma}
\label{ratio}
Given positive numbers $a_1,\dots,a_n$ and $b_1,\dots,b_n$ we have 
$$\frac{a_1+\cdots+a_n}{b_1+\cdots+b_n}\leq\max\left\{\frac{a_1}{b_1},\dots,
\frac{a_n}{b_n}\right\}.$$
\end{lemma} 

\begin{theorem}
\label{est}
For a $d$-polytope $\Delta$ simple in edges and for $1\leq k<l\leq d/2$
there is the following upper bound for the average number of $k$-dimensional
faces lying on a $l$-dimensional face of $\Delta$: 
$${n-k\choose n-l}\frac{{[d/2]\choose k}+{[(d+1)/2]\choose k}}
{{[d/2]\choose l}+{[(d+1)/2]\choose l}}$$
\end{theorem}

\proof First note that for any $k$-dimensional face $F$ of $\Delta$
($k\geq 1$) there are exactly ${n-k\choose n-l}$ faces of dimension $l$
containing $F$. Therefore it suffices to estimate the ratio $f_k/f_l$. Using
theorem \ref{ineq3} and the relation between the $f$- and $h$-vectors we get
$$f_k\leq\sum_{m\geq d/2}h_m\left[{m\choose k}+{d-m\choose k}\right].$$
$$f_l\leq\sum_{m\geq d/2}h_m\left[{m\choose l}+{d-m\choose l}\right].$$  
So we can use lemma \ref{ratio} to estimate $f_k/f_l$ from above. We
finally get
$$\frac{f_k}{f_l}\leq\frac{{[d/2]\choose k}+{[(d+1)/2]\choose k}}
{{[d/2]\choose l}+{[(d+1)/2]\choose l}}.\ \Box$$ 

\section{Polytopes with infrequent singularities}

\paragraph{Single nonsimple vertex.}
Consider a polytope $\Delta$ with only one nonsimple vertex $v$. Let
$\Sigma$ be a simple polytope obtained from $\Delta$ by cutting off the
vertex $v$. Denote the only inserted face of $\Sigma$ by $\Gamma$.

\begin{lemma}
\label{ins}
Suppose that an element $\alpha\in A^k(\Sigma)$ is such that
$\d_{\Gamma'}\alpha=0$ for each facet $\Gamma'$ of $\Sigma$ that
does not intersect $\Gamma$. Then $\alpha$ is divisible by $\d_\Gamma$
in $A(\Sigma)$. 
\end{lemma}

\proof
Introduce a general linear function $l$ on $\Sigma$ such that all vertices
of $\Gamma$ are lower than all others. Consider the separatrix decomposition
of $\alpha$ with respect to $l$. Let $F$ be the highest separatrix.
Assume that $F$ does not belong to $\Gamma$. Then there exists a facet
$\Gamma'$ of $\Sigma$ such that $\Gamma'$ passes through the top vertex of
$F$, does not contain $F$ and does not intersect $\Gamma$.
We know that $\d_{\Gamma'}\alpha=0$ (or, equivalently,
$\alpha_{(\Gamma')}=0$). On the other hand, it is easy to see
that $(\d_F)_{(\Gamma')}$ is the highest separatrix operator
of $\alpha_{(\Gamma')}$ with respect to the general linear function
$l|_{\Gamma'}$ on $\Gamma'$. Contradiction.

So any face $F$ from the separatrix decomposition of $\alpha$ belongs to
$\Gamma$. Hence $\alpha=\d_\Gamma\beta$ for some $\beta$. $\Box$

\begin{theorem}
\label{one-s}
Suppose $\alpha L_\Delta^{d-2k}=0$ where $\alpha\in A^k(\Sigma)$. 
Then $\alpha$ is divisible by $\d_\Gamma$ in $A(\Sigma)$.
\end{theorem}

\proof Place the origin at $v$. Then the number $H_{\Gamma'}(\Delta)$ is
nonzero only if $\Gamma'$ does not intersect $\Gamma$. Take such a facet
$\Gamma'$. Project the equality $\alpha L_\Delta^{d-2k}(\Vol_\Sigma)=0$
to the facet $\Gamma'$. We will get
$\alpha_{(\Gamma')}L_{\Gamma'}^{d-2k}(\Vol_{\Gamma'})=0$
(we know that $(L_\Delta)_{(\Gamma')}=L_{\Gamma'}$). But 
this is the primitivity condition with respect to $\Gamma'$. By theorem
\ref{HR}, we have 
$$(-1)^k\alpha_{(\Gamma')}^2L_{\Gamma'}^{d-1-2k}(\Vol_{\Gamma'})\geq 0.$$ 
Multiply this inequality by $H_{\Gamma'}(\Delta)>0$ and sum up over all
facets $\Gamma'$ not intersecting $\Gamma$. We get
$(-1)^k\alpha^2L_\Delta^{d-2k}(\Vol_\Sigma)\geq 0$.
But this is an equality according to our assumption on $\alpha$.
Therefore for all $\Gamma'$ (such that $\Gamma'\cap\Gamma=\varnothing$) we
have $(-1)^k\alpha_{(\Gamma')}^2L_{\Gamma'}^{d-1-2k}(\Vol_{\Gamma'})=0$. 
Since $\alpha_{(\Gamma')}$ is $\Gamma'$-primitive,
$\alpha_{(\Gamma')}=0$ in $A^k(\Gamma')$ or, equivalently,
$\d_{\Gamma'}\alpha=0$. By lemma \ref{ins}, $\alpha$ is divisible by
$\d_\Gamma$. $\Box$

The method we used in this proof is very similar to those of Aleksandrov
\cite{Al} and McMullen \cite{McM}. 

\paragraph{Infrequent singularities.}
We need the following simple fact:

\begin{lemma}
\label{reldiff}
Let $P$ and $Q$ be homogeneous polynomials on the same vector space.
A differential operator $\beta$ with constant coefficients such that
$P=\beta Q$ exists if and only if $\alpha Q=0$ implies $\alpha P=0$ for each
differential operator $\alpha$ with constant coefficients. 
\end{lemma}
       
\proof The part ``only if'' is obvious.
Now assume that from $\alpha Q=0$ it always follows that $\alpha P=0$.
Any operator of order $>\deg(Q)$ annihilates $P$ so $\deg(P)\leq\deg(Q)$.
Denote by $A$ the quotient of the polynomial algebra with respect to the
ideal annihilating $Q$. Let $W$ be the hyperplane in $A^{\deg(P)}$
consisting of all operators $\alpha$ such that $\alpha P=0$. Denote by
$\beta\in A^{\deg(Q)-\deg(P)}$ any generator of the one-dimensional
orthogonal complement to $W$ with respect to the nondegenerate pairing
$(\alpha,\beta)\mapsto \alpha\beta Q$.

The polynomial $P$ can be viewed as a linear functional on $A^{\deg(P)}$.
The functionals $P$ and $\beta Q$ have the same zero level. Therefore
they are proportional. $\Box$

Let $\Sigma$ be a simple polytope. Suppose we want to prove that a
polynomial $P$ in support numbers of $\Sigma$ has the form
$\beta\Vol$ for some $\beta\in A(\Sigma)$. Then by lemma \ref{reldiff}
it is enough to verify that the ideal $J=\{\alpha\in\Diff|\ \alpha\Vol=0\}$
annihilates $P$. It suffices to check that all the generators of $J$
send $P$ to zero. Namely, for each collection $\Gamma_1,\dots,\Gamma_k$
of facets of $\Sigma$ with empty intersection we should show that
$\d_{\Gamma_1}\cdots\d_{\Gamma_k}P=0$, and for each point $a\in\R^d$
we need to prove that $L_a P=\sum H_\Gamma(a)\d_\Gamma(P)=0$. 

Let $\Delta$ be a polytope simple in edges such that no its facet contains
more than one nonsimple vertex. We call such a polytope a {\em polytope with
infrequent singularities}. Denote by $\Sigma$ the standard resolution of
$\Delta$. Define the space $I^k$ as the subspace of $A^k(\Sigma)$ generated by
all elements of the form $\d_\Gamma\alpha$ where $\Gamma$ is an inserted
facet of $\Sigma$ and $\alpha\in A^{k-1}(\Sigma)$.

\begin{theorem}
\label{ker}
If $\Delta$ is a polytope with infrequent singularities, then the operator 
(of multiplication by) $L_\Delta:A^k(\Sigma)\to A^{k+1}(\Sigma)$ has the
kernel $I^k$ for $k<(d-1)/2$. 
\end{theorem}

\proof
Let $\Gamma$ be an inserted facet. Then $L_\Delta\d_\Gamma=0$ since
$(L_\Delta)_{(\Gamma)}=L_v=0$. Hence the subspace $I^k$ lies in the
kernel of $L_\Delta$. 

Let us prove the opposite inclusion. 
We will carry on the induction on the number of nonsimple vertices of
$\Delta$. If there is only one nonsimple vertex, then the theorem follows
from theorem \ref{one-s}.
Now let $v$ be an arbitrary nonsimple vertex of $\Delta$. Cut it off. 
We get another polytope $\Theta$ with infrequent singularities. Denote by 
$\Gamma$ the inserted facet of $\Theta$ corresponding to $v$. We can assume
that $\Gamma$ coincides with the corresponding inserted facet of $\Sigma$. 
In particular, $\Sigma$ is a standard resolution of $\Theta$.

Let $\alpha\in A^k(\Sigma)$ satisfy $L_\Delta\alpha=0$. Note that
$L_\Delta-L_{\Theta}=c\d_{\Gamma}$ where $c$ is a positive number. Thus
$L_{\Theta}\alpha=-c\d_\Gamma\alpha$. 

From the equality $L_\Delta\alpha=0$ we get an analogous relation
$(L_\Delta)_{(\Gamma')}\alpha_{(\Gamma')}=0$ for each non-inserted facet
$\Gamma'$ of $\Sigma$. The operator $(L_\Delta)_{(\Gamma')}$ corresponds to
a facet of $\Delta$ with at most one nonsimple vertex and the standard
resolution $\Gamma'$. By theorem \ref{one-s} for any facet $\Gamma'$
intersecting $\Gamma$ the element $\alpha_{(\Gamma')}$ is divisible
by $(\d_\Gamma)_{(\Gamma')}$, i.e. $\d_{\Gamma'}\alpha$ is divisible by
$\d_\Gamma$. Set $P_{\Gamma'}=\d_{\Gamma'}\alpha\Vol$ if $\Gamma'$
intersects $\Gamma$ and $P_{\Gamma'}=0$ otherwise.

It is easy to verify that the polynomials $P_{\Gamma'}$ are related
as follows: $\d_{\Gamma'}P_{\Gamma''}=\d_{\Gamma''}P_{\Gamma'}$.
Therefore there exists a polynomial $P$ such that $P_{\Gamma'}=\d_{\Gamma'}P$
for all facets $\Gamma'$. We want to prove that $P=\beta\Vol$ for
some differential operator $\beta$ with constant coefficients. By lemma
\ref{reldiff} it is enough to show that
\begin{itemize}
\item $\d_{\Gamma'_1}\cdots\d_{\Gamma'_k}P=0$ if
$\Gamma'_1\cap\cdots\cap\Gamma'_k=\varnothing$,
\item $\sum H_{\Gamma'}(a)\d_{\Gamma'}P=0$ for each $a\in\R^d$.
\end{itemize}
The first condition is obvious. The second condition follows from the
equation $\sum H_{\Gamma'}(a)\d_{\Gamma'}\alpha=0$ in $A(\Sigma)$.
Note that by theorem \ref{one-s} each summand $\d_{\Gamma'}\alpha$ is
divisible by an operator of some inserted facet. 
To obtain the second condition on $P$ equate to zero the terms with
$\d_\Gamma$ only (using lemma \ref{indep}).

Thus we have $P=\beta\Vol$. But $\d_{\Gamma'}P=0$ for any $\Gamma'$
not intersecting $\Gamma$. From lemma \ref{ins}
it follows that $\beta$ is divisible by $\d_\Gamma$ and, in particular,
$L_\Delta\beta=0$. By definition of $\beta$
each derivative of $\alpha_{(\Gamma)}\Vol_\Gamma$ coincides with the
corresponding derivative of $\beta_{(\Gamma)}\Vol_\Gamma$.
Therefore $\d_\Gamma\alpha=\d_\Gamma\beta$ and $L_\Theta\alpha=L_\Theta\beta$. 

Now let $\gamma=\alpha-\beta$. We know that $L_\Theta\gamma=0$. By the
induction hypothesis $\gamma\in I^k$. Thus $\alpha=\beta+\gamma\in I^k$.
$\Box$

\section{Consequences for combinatorial intersection cohomology}

In this section, we will give an interpretation of theorem \ref{ker}
in terms of the combinatorial intersection cohomology.
First let us recall briefly some basic definitions.

\paragraph{Fans.} To each face $F$ of a $d$-polytope associate the
{\em normal cone} $C_F\subset\R^{d*}$ consisting of linear functionals
on $\R^d$ that achieve their maximal values somewhere on $F$. The set
of normal cones to all the faces of a polytope $\Delta$ is called the {\em
dual fan} of $\Delta$.

{\em A fan} in a real vector space $V$ is a collection $\Phi$
of convex polyhedral cones with vertex at the origin such that 
\begin{itemize}
\item for every cone $\sigma\in\Phi$ all the faces of $\sigma$ belong to
$\Phi$,
\item the intersection of two cones in $\Phi$ is their common face. 
\end{itemize}
A fan is said to be {\em simplicial} if all its cones are simplicial.
A fan is {\em complete} if the union of all its cones is the whole space $V$.

The dual fan of a polytope in $\R^d$ is a complete fan in $\R^{d*}$.
It is simplicial if and only if the corresponding polytope is simple.

\paragraph{Toric varieties.} Fix a lattice $\Omega$ in a vector space $V$.
A fan in $V$ is said to be {\em rational} if all its rays
(i.e., one-dimensional cones) are spanned by
lattice vectors. For each rational fan $\Phi$ one defines the corresponding
{\em toric variety} $X$. This is a complex algebraic variety
with an algebraic action of the complex torus $\T=(V\otimes\C)/i\Omega$.
If $\Phi$ is complete, then $X$ is compact; if $\Phi$ is simplicial, then
$X$ is an orbifold; if $\Phi$ is dual to a polytope, then $X$ is
projective.

Suppose that the dual fan $\Phi$ of a polytope $\Delta$ is rational. Then
the intersection cohomology Betti numbers of the corresponding toric
variety $X$ are combinatorial
invariants of $\Delta$. The intersection cohomology of $X$
can be described explicitly in terms of $\Phi$ only \cite{BL,BBFK}.
This description makes sense even then $\Phi$ is nonrational and
there is no corresponding toric variety. 

\paragraph{Combinatorial intersection cohomology.} Following \cite{BL,BBFK}
we will define the (combinatorial) intersection cohomology of a fan.
A fan $\Phi$ can be considered as a finite topological space whose
open subsets are subfans. Every cone $\sigma\in\Phi$ has a unique
minimal neighborhood $[\sigma]$ consisting of $\sigma$ and all its faces.

Let us define a sheaf of rings $\Oc_\Phi$ on $\Phi$. Sections of
$\Oc_\Phi$ over a subfan $\Ups$ are continuous functions on $\bigcup\Ups$
that are polynomial on each cone of $\Ups$. It is not hard to verify that
$\Oc_\Phi$ is flabby if and only if $\Phi$ is simplicial.

A graded sheaf $\Mc_\Phi$ of $\Oc_\Phi$-modules is called {\em basic} if it
satisfies the following conditions:
\begin{itemize}
\item {\em Normalization}: $\Mc_\Phi([0])=\R$.
\item {\em Pointwise freeness}: $\Mc_\Phi[\sigma]$ is a
free $\Oc_\Phi[\sigma]$-module for any $\sigma\in\Phi$.
\item {\em Flabbyness}: the sheaf $\Mc_\Phi$ is flabby. For that it is
enough to require that for any cone $\sigma\in\Phi$ the restriction
map $\Mc_\Phi[\sigma]\to\Mc_\Phi(\d\sigma)$ be surjective. 
\item {\em Minimality}: the module $\Mc_\Phi[\sigma]$ is a minimal
free $\Oc_\Phi[\sigma]$-module satisfying the previous condition.
\end{itemize}
It is clear that a basic sheaf exists and is unique up to isomorphism.  
The space of global sections $M_\Phi=\Gamma(\Phi,\Mc_\Phi)$ is a
combinatorial analog of the equivariant intersection cohomology of a toric
variety.

Let $O_\Phi=\Gamma(\Phi,\Oc_\Phi)$ be the space of global
continuous piecewise polynomial functions on $\Phi$. The space
$M_\Phi$ is a $O_\Phi$-module. Consider the ideal $O_\Phi^+$ in $O_\Phi$
generated by all global linear functions.
Denote the quotient module $M_\Phi/O_\Phi^+M_\Phi$ by $\overline{M}_\Phi$.
This is the {\em intersection cohomology} of $\Phi$. 

\paragraph{An analog of the Hard Lefschetz theorem.} Let $\Phi$ be the dual
fan of a polytope $\Delta\subset\R^d$. 
For any linear functional $\xi\in\R^{d*}$ denote by $S_\Delta(\xi)$ the
maximum of $\xi$ restricted to $\Delta$. The function $S_\Delta$
is piecewice linear with respect to $\Phi$. Hence it lies in $O_\Phi$.
The following analog of the Hard Lefschetz theorem holds
for integral polytopes \cite{BL,BBFK} and is believed to be true for
arbitrary polytopes.

\begin{conjecture}
\label{HL}
For $k<d/2$ the operator of multiplication by $S_\Delta^{d-2k}$ establishes
an isomorphism between $\overline{M}^k_\Phi$ and $\overline{M}^{d-k}_\Phi$. 
In particular, the multiplication by $S_\Delta$ is an embedding of
$\overline{M}^k_\Phi$ to $\overline{M}^{k+1}_\Phi$. 
\end{conjecture}

Denote the combinatorial intersection Betti numbers $\dim\overline{M}^k_\Phi$
by $Ih_k(\Delta)$. If conjecture \ref{HL} is true, then
$$Ih_0(\Delta)\leq Ih_1(\Delta)\leq\dots\leq Ih_{[d/2]}(\Delta).$$
It is proven in \cite{BL,BBFK} that $Ih_k(\Delta)=Ih_{d-k}(\Delta)$
(an analog of Poincar\'e duality).

\paragraph{Cohomology of a simplicial fan.} Let $\Psi$ be a simplicial
fan. Then the basic sheaf $\Mc_\Psi$ coincides with $\Oc_\Psi$.
Therefore $M_\Psi=O_\Psi$ is the space of all piecewise polynomial
functions on $\Psi$. Now assume that $\Psi$ is dual to a simple polytope
$\Sigma$. 

\begin{proposition}
\label{iso}
There is a natural isomorphism between $\overline{O}_\Psi=O_\Psi/O^+_\Psi$
and $A(\Sigma)$. This isomorphism takes $S_\Sigma$ to $L_\Sigma$. 
\end{proposition}

\proof For a ray $\rho\in\Psi$ denote by
$\chi^\rho$ a piecewise linear function that is zero on all the rays of
$\Psi$ but $\rho$. The function $\chi^\rho$ will be called a {\em
characteristic function} of $\rho$. Characteristic function of $\rho$ is
unique up to constant factor (of course, we assume that $\chi^\rho\ne 0$).
Note that if rays
$\rho_1,\dots,\rho_k$ do not lie in a common cone of $\Psi$, then
$\chi^{\rho_1}\cdots\chi^{\rho_k}=0$.

Now pass to $A(\Sigma)$. This algebra is generated by differentiations
$\d_\Gamma$ where $\Gamma$ are facets of $\Sigma$. The relations in
$A(\Sigma)$ are generated by the following two groups \cite{Tim}:
\begin{itemize}
\item {\em Incidence relations:}
if $\Gamma_1\cap\cdots\cap\Gamma_k=\varnothing$, then
$\d_{\Gamma_1}\dots\d_{\Gamma_k}=0$,
\item {\em Translation invariance:}
if $a\in\R^d$ is a point, then $\sum H_\Gamma(a)\d_\Gamma=0$. 
\end{itemize}
The incidence relations provide the homomorphism $\phi:O_\Psi\to
A(\Sigma)$ that takes $\chi^\rho$ to $\d_\Gamma$ where $\rho$ is the normal
cone of $\Gamma$.
This homomorphism is clearly onto. Translation invariance relations
determine the kernel of $\phi$. It is generated by all points $a\in\R^d$
considered as linear functions on $\R^{d*}$. But this is the same as
$O^+_\Psi$. Therefore $\overline{O}_\Psi$ and $A(\Sigma)$ are isomorphic.
It is easy to see that the isomorphism thus constructed takes $S_\Sigma$
to $L_\Sigma$. $\Box$

\paragraph{Cohomology of a single cone.} Consider a noncomplete fan
$[\sigma]$ that consists of a $d$-dimensinal cone $\sigma$ and all its faces.
We will assume that $\sigma$ contains no nontrivial vector subspace.
The fan $[\sigma]$ corresponds (in the rational case) to an affine toric
variety. The ring $O_{[\sigma]}$ consists of all polynomials on $\R^{d*}$.

Let us study the restriction $M_{[\sigma]}\to M_{\d\sigma}$.
One can project the fan $\d\sigma$ to a complete fan $\overline\sigma$ in
a proper hyperplane passing through the origin. Clearly, the fan 
$\overline\sigma$ is the dual fan of a polytope $\Lambda$.
If $\sigma$ is the cone over a polytope $\Lambda^*$, then $\Lambda$
is dual to $\Lambda^*$ in combinatorial sense, i.e., there is
a one-to-one correspondence between proper faces of $\Lambda$ and
$\Lambda^*$ reversing inclusions. 
The restriction
$\Mc|_{\d\sigma}$ defines a basic sheaf on the fan $\overline\sigma$
and we have $M_{\overline\sigma}=M_{\d\sigma}$.
From the minimality condition for the sheaf $\Mc_{[\sigma]}$ and Nakayama's
lemma it follows that the map
$$\overline{M}_{[\sigma]}=M_{[\sigma]}/O_{[\sigma]}^+M_{[\sigma]}\to
M_{\d\sigma}/O_{[\sigma]}^+M_{\d\sigma}$$
is an isomorphism of vector spaces over the field
$\overline{O}_{[\sigma]}=O_{[\sigma]}/O^+_{[\sigma]}=\R$.

Let us introduce coordinates $(x_1,\dots,x_d)$ in the space $\R^{d*}$
so that the fan $\overline\sigma$ lies in the coordinate hyperplane $x_d=0$.
Then the linear function $x_d$ restricted to $\d\sigma$ and projected
to $\overline\sigma$ acts on the fan $\overline\sigma$ as the Lefschetz
operator of $\Lambda$.

\begin{proposition}
\label{kerres}
Suppose the fan $\overline\sigma$ satisfies conjecture \ref{HL}
(with respect to the function $S_\Lambda$).
Then the kernel of the restriction map $M_{[\sigma]}\to M_{\overline\sigma}$
contains no elements of degree $\leq d/2$.
\end{proposition}

\proof
Consider the ring $o=\R[x_1,\dots,x_{d-1}]$ and the maximal ideal
$o^+$ in it generated by all the degree-one elements (i.e., by all
linear functions). The cohomology
space of the fan $\overline{\sigma}$ is $\overline{M}_{\overline\sigma}=
M_{\overline\sigma}/o^+M_{\overline\sigma}$.

Let $\phi_1,\dots,\phi_s$ be free generators of the module 
$M_{[\sigma]}$ over the ring $O_{[\sigma]}=\R[x_1,\dots,x_d]$. Denote by
$\bar\phi_i$ the image of the element $\phi_i$ under the restriction
$M_{[\sigma]}\to M_{\d\sigma}$. The elements $\bar\phi_i$ obviously
give rise to a basis in the vector space
$M_{\d\sigma}/O^+_{[\sigma]}M_{\d\sigma}\cong\overline{M}_{[\sigma]}$.

Now assume that an element $\phi=a_1\phi_1+\cdots+a_s\phi_s$ homogeneous
of degree $k\leq d/2$ restricts to zero, i.e.,
$a_1\bar\phi_1+\cdots+a_s\bar\phi_s=0$ in $M_{\overline{\sigma}}$. Reduce
the latter relation modulo the ideal $o^+$. After that all the coefficients
$a_i$ become polynomials in $x_d$. Let $x_d^t$ be the least power of $x_d$
that divides all the polynomials $a_i$. As $k\leq d/2$, the operator
of multiplication by $x_d^t$ is an injective map from
$\overline{M}_{\overline\sigma}^{k-t}$ to $\overline{M}_{\overline\sigma}^k$.
This follows from conjecture \ref{HL} for $\overline\sigma$.
Hence it is possible to divide our relation by $x_d^t$.
We get $b_1\bar\phi_1+\cdots+b_s\bar\phi_s=0$ where at least 
one coefficient $b_i$ is not divisible by $x_d$. But then we reduce 
this relation modulo $(x_d)$ to obtain a nontrivial linear combination 
of the elements $\bar\phi_i$ in
$M_{\d\sigma}/O^+_{[\sigma]}M_{\d\sigma}=\overline{M}_{[\sigma]}$ that equals
to zero. This contradicts the statement that (the classes of) the elements
$\bar\phi_i$ constitute a basis of the vector space
$\overline{M}_{[\sigma]}$. $\Box$

Let us point out just another corollary from the Hard Lefschetz theorem
for the fan $\overline\sigma$:

\begin{proposition}
Suppose the fan $\overline\sigma$ satisfies conjecture \ref{HL}.
Then $\overline{M}_{[\sigma]}^k=0$ for $k>d/2$.
\end{proposition}

\proof
Indeed, $\overline{M}_{[\sigma]}=\overline{M}_{\overline\sigma}/(x_d)
\overline{M}_{\overline\sigma}$.
From the Hard Lefschetz theorem it follows that all the elements of degree
$k>d/2$ in the space $\overline{M}_{\overline\sigma}$ lie in the image of the Lefschetz
operator, i.e., in $(x_d)\overline{M}_{\overline\sigma}$. $\Box$

The same arguments help to compute dimensions of the homogeneous 
components $\overline{M}_{[\sigma]}^k$ ($k\leq d/2$). From the Lefschetz
decomposition for the polytope $\Lambda$ we obtain
$$\dim(\overline{M}^k_{[\sigma]})=
Ih_k(\Lambda)-Ih_{k-1}(\Lambda),\quad k\leq d/2.$$
We see that dimension of the vector space $\overline{M}_{[\sigma]}^k$ is a
combinatorial invariant of the polytope $\Lambda$. Denote this dimension by
$Ig_k(\Lambda)$. Thus
$$Ig_k(\Lambda)=\left\{\begin{array}{cl}
Ih_k(\Lambda)-Ih_{k-1}(\Lambda),& k\leq d/2,\\
0,& k>d/2.\end{array}\right.$$

\paragraph{Links and the computaion of $Ih_k(\Delta)$.}
Let $F$ be a face of a polytope $\Delta$. Denote by $N(F)$ the
orthogonal complement to the plane of the face $F$ passing through an
interior point of $F$. In a neighborhood of this point the intersection 
$N(F)\cap\Delta$ looks like a cone over a polytope $\Lambda(F)$.
The polytope $\Lambda(F)$ is called the {\em link} of the face $F$. 
Let $\sigma\in\Phi$ be the normal cone of the face $F$. Then the 
dual fan of the polytope $\Lambda(F)$ is $\overline\sigma$.

Consider the generating function $IH_\Delta(t)=\sum Ih_k(\Delta)t^k$.
For a polytope $\Lambda$ denote $IG_\Lambda(t)=\sum Ig_k(\Lambda)t^k$.
In \cite{BL} it is proved that the Hard Lefschetz theorem \ref{HL} for all
links
of $\Delta$ would imply the following formula for the cohomology of $\Delta$:
$$IH_\Delta(t)=\sum_F (t-1)^{\dim F} IG_{\Lambda(F)}(t)\eqno{(*)}$$
($F$ runs over all the faces of $\Delta$).
The proof of this formula splits naturally into two parts. The first
part is the computation of the cohomology for the fan $[\sigma]$. 
This computation is already done, it relies on the Hard Lefschetz 
theorem. The second part reduces the global cohomology to the local 
cohomology, i.e., to the cohomology of the fans $[\sigma]$. This part
does not depend on the Hard Lefschetz theorem.

All links of a polytope simple in edges are simple. Therefore for
a polytope $\Delta$ simple in edges the above formula for $IH_\Delta$ is
true. It is easy to verify using this formula that for $k>d/2$ we have
$Ih_k(\Delta)=h_k(\Delta)$.

\paragraph{Stanley's generalized $h$-vector.} For an arbitrary polytope
$\Delta$ Stanley in \cite{St2} defined the {\em generalized
$h$-vector} $Gh$ by the combinatorial recurrent formula $(*)$ (where
$Ih$-vector is replaced of course by $Gh$). As we
saw, for rational polytopes, simple polytopes and polytopes simple in
edges we have $Gh=Ih$. It is believed that this holds in general.

Stanley proved that $Gh_k(\Delta)=Gh_{d-k}(\Delta)$ and conjectured that
$$Gh_0\leq Gh_1\leq\cdots\leq Gh_{[d/2]}.$$
We will see that this is true for polytopes with infrequent singularities.

\paragraph{Cohomology of polytopes simple in edges.} Let $\Delta$ be
a $d$-polytope simple in edges. The dual fan $\Phi$ of $\Delta$ has the
following
property: its $(d-1)$-skeleton $\Ups$ is simplicial. Now consider a
standard resolution $\Sigma$ of $\Delta$. Its dual fan $\Psi$ is a
simplicial subdivision of $\Phi$. It an easy exercise to describe
this subdivision explicitly. Note that $\Ups$ is a subfan both
in $\Phi$ and $\Psi$. 

We know already how the cohomology $\overline{M}_\Psi=\overline{O}_\Psi$
looks like. Now we need the following theorem proved in \cite{BL,BBFK}:

\begin{theorem}
Let $j:\Psi\to\Phi$ be the natural map that takes each cone $\sigma\in\Psi$
to the minimal cone of $\Phi$ where $\sigma$ lies. There exists a
(noncanonical) embedding $\Mc_\Phi\to j_*\Oc_\Psi$ that preserves the
structure of $\Oc_\Phi$-modules.
\end{theorem}
                        
Fix any such embedding. Then $\overline{M}_\Phi$ can be viewed as a subspace
of $\overline{O}_\Psi=A(\Sigma)$. Denote by $\overline{N}$ the kernel of the
restriction homomorphism $\bar R:\overline{O}_\Psi\to\overline{O}_\Ups$. 

\begin{proposition}
\label{comM}
For $k\leq d/2$ there is the following decomposition
$$\overline{O}_\Psi=\overline{M}_\Phi^k\oplus\overline{N}^k.$$
\end{proposition}

\proof
It is enough to prove that on $\overline{M}_\Phi^k$ the operator $\bar R$ is an
isomorphism. Moreover, it is sufficient to prove only injectiveness
(surjectiveness is clear).

Take an element $\bar\phi\in\Ker(\bar R)$ of degree $k$ and consider its
representative $\phi\in M^k_\Phi$. Then $\phi$ restricts to $O^+_\Ups$, i.e.,
$\phi=l_1\theta_1+\cdots+l_r\theta_r$ on $\Ups$ where $l_i$ are linear
functions on $\R^{d*}$ and $\theta_i\in O^{k-1}_\Ups$. Since $\Mc_\Phi$ is
flabby each $\theta_i$ comes from an element $\phi_i\in M^{k-1}_\Phi$. The
function $\phi'=\phi-l_1\phi_1-\cdots-l_r\phi_r\in M^k_\Phi$ represents
$\bar\phi$ and is zero on $\Ups$.

Let $\sigma$ be a cone of dimension $d$ in $\Phi$ such that $\phi'$ is
nonzero on $\sigma$. The restriction of $\phi'$ to the subfan $[\sigma]$ is a 
nonzero element of the kernel of the restriction map $M_{[\sigma]}\to
M_{\d\sigma}$. But according to proposition \ref{kerres} the degree
of such an element can not be $\leq d/2$. Contradiction. $\Box$

\begin{proposition}
\label{transl}
Under the identification $\overline{O}_\Psi=A(\Sigma)$ the subspace
$\overline{N}^k$ coincides with $I^k$. 
\end{proposition}

\proof We know that $\d_\Gamma$ corresponds to $\chi^\rho$ where $\rho$
is the dual ray of $\Gamma$. Thus it is enough to show that $\overline{N}$
is generated (as a $\overline{O}_\Psi$-module) by the characteristic functions
of rays $\rho\in\Psi-\Phi$. It is clear that all such functions belong
to $\overline{N}$.

Now take $\bar\phi\in\overline{N}$ and its representative $\phi\in O_\Psi$
such that $\phi=0$ on $\Ups$ (such a representative obviously exists).
Let us define an {\em inserted cone} as a cone from $\Psi-\Phi$.
Inserted rays correspond to inserted facets of $\Sigma$. 
Suppose $\phi$ is nonzero on an inserted ray $\rho$. Then we can
subtract from $\phi$ an appropriate multiple of $\chi^\rho$
(i.e., a function of the form $\psi\chi^\rho$, $\psi\in O_\Psi$) so that
the result becomes zero on $\rho$ and remains the same on all other
rays of $\Psi$. Thus we can reduce $\phi$ (modulo characteristic functions
of inserted rays) to a function $\phi'$ that is zero on all the rays
of $\Psi$. Now repeat this procedure with 2-dimensional cones.
Suppose $\phi'$ is nonzero on an inserted 2-dimensional cone $\tau$
bounded by rays $\rho_1$ and $\rho_2$. One of these rays is inserted.
Subtract from $\phi'$ an appropriate multiple of
$\chi^{\rho_1}\chi^{\rho_2}$ to obtain a function that is zero on $\tau$
and the same as $\phi'$ on all other 2-dimensional cones of $\Psi$.
Continueing this process we reduce $\phi$ to zero modulo characteristic
functions of inserted rays. $\Box$

\paragraph{Hard Lefschetz for polytopes with infrequent singularities}
Combining some previous results (theorem \ref{ker}, propositions
\ref{comM}, \ref{transl}) we obtain the following theorem:

\begin{theorem}
Let $\Delta$ be a polytope with infrequent singularities and $\Phi$ the
dual fan of it. The multiplication by $S_\Delta$ establishes an
embedding of $\overline{M}^k_\Phi$ to $\overline{M}^{k+1}_\Phi$ for
$k<(d-1)/2$. In particular, 
$$Gh_0(\Delta)\leq Gh_1(\Delta)\leq\cdots\leq Gh_{[d/2]}(\Delta),$$
$$h_{[d/2]}(\Delta)\geq h_{[d/2]+1}(\Delta)\geq\cdots\geq h_d(\Delta).$$
\end{theorem}

\end{document}